\documentclass[12pt]{article}
\usepackage{amsmath}
\usepackage{color}
\usepackage{fancyhdr}
\usepackage{verbatim}
\usepackage{amsfonts,longtable,mathrsfs}
\usepackage{epsfig}
\usepackage{amsfonts}
\usepackage{color}
\usepackage[numbers,sort&compress]{natbib}

\def\R{\hbox{{\rm I}\kern-0.2em{\rm R}\kern0.2em}}

\def\bn{\begin{equation}}
\def\en{\end{equation}}
\def\bny{\begin{eqnarray}}
\def\eny{\end{eqnarray}}
\def\be{\begin{eqnarray*}}
\def\ee{\end{eqnarray*}}
\def\bc{\begin{center}}
\def\ec{\end{center}}
\def\X{{\cal X}} 

\def\({\left(}
\def\){\right  )}
\def\[{\left[}
\def\]{\right]}
\def\bc{\begin{center}}
\def\ec{\end{center}}

\newtheorem{dfn}{Definition}[section]
\newtheorem{thm}{Theorem}[section]
\newtheorem{rem}{Remark}[section]
\newtheorem{pro}{Proposition}[section]

\newtheorem{cor}{Corollary}[section]
\newtheorem{lem}{Lemma}[section]
\newtheorem{exm}{Example}[section]

\def\bn{\begin{equation}}
\def\en{\end{equation}}
\def\bny{\begin{eqnarray}}
\def\eny{\end{eqnarray}}
\def\be{\begin{eqnarray*}}
\def\ee{\end{eqnarray*}}
\def\bdn{\begin{dfn}}
\def\edn{\end{dfn}}
\def\btm{\begin{thm}}
\def\etm{\end{thm}}
\def\bpf{\begin{proof}}
\def\epf{\end{proof}}
\def\bpn{\begin{pro}}
\def\epn{\end{pro}}
\def\brk{\begin{rem}}
\def\erk{\end{rem}}
\def\bcy{\begin{cor}}
\def\ecy{\end{cor}}
\def\blm{\begin{lem}}\def\elm{\end{lem}}
\def\bex{\begin{exm}}
\def\eex{\end{exm}}
\def\X{{\cal X}}  
  
 \def\R{{\hat R}}
\begin{document}

\bc {\bf The Invariance and Conservation Laws of fourth-order Difference Equations
  }\ec
\medskip
\bc
M. Folly-Gbetoula$^{1,}$\footnote{Mensah.Folly-Gbetoula@wits.ac.za} and  A. H. Kara$^{1,2,}$\footnote{Abdul.Kara@wits.ac.za }
\\$^1$School of Mathematics, University of the Witwatersrand, Johannesburg, South Africa.\\
$^2$ Department of Mathematics and Statistics, King Fahd University of Petroleum and Minerals, Dhahran, Saudi Arabia.
\ec
\begin{abstract}
We consider difference equations of order four and determine the one parameter Lie group of transformations (Lie symmetries) that leave them invariant. We introduce a technique for finding their first integrals and discuss the association between the symmetries and first integrals as well as the notion of multipliers (related to conservation laws) for difference equations.
\end{abstract}
\textbf{Key words}: Difference equations; symmetries; reduction; group invariant solutions
\section{Introduction} \setcounter{equation}{0}
Difference equations ($\Delta$Es), also known as recurrence relations, are mathematical models of evolution of discrete phenomena. They have been studied by a number of authors \cite{Hydon1, l1, s, FK1, FK2} and it is now known that they can be analysed using the symmetry methods as was done for differential equations. The symmetry method for differential equations is well-documented and can be found in \cite{9,10} among others. The main ideas and methods for symmetry analysis of difference equations were introduced by Hydon in \cite{11} and Winternitz et al.
in \cite{W}.
\subsection{Overview about Lie analysis of $\Delta E$'s}
Let us consider an $N$th-order difference equation
\begin{equation}\label{general}
u_{n+N}=\omega(n, u_n, u_{n+1}, \dots, u_{n+N-1})
\end{equation}
for some function $f$ (we will assume that $\partial \omega/ \partial {u_n}\neq 0)$ and the point transformations
\begin{equation}\label{Gtransfo}
\Gamma _{\epsilon}: \textbf{x} \mapsto \tilde{\textbf{x}}(\textbf{x};\epsilon),
\end{equation}
where $\textbf{x}= (x_1,x_2,\dots,x_p)$ are the continuous variables.
$\Gamma$ is a one-parameter Lie group of transformations if:
\begin{itemize}
\item
$\Gamma _0$ is the identity map, so that $\tilde{\textbf{x}}=\textbf{x}$ when $\epsilon=0$
\item
$\Gamma _{a}\Gamma _{b} = \Gamma _{a+b}$ for every $a$ and $b$ sufficiently close to $0$
\item
Each $\tilde{x}_i$ can be represented as a Taylor series in $\epsilon$ (in a neighborhood of $\epsilon =0$ that is determined by $\textbf{x}$), and therefore
\begin{equation}\label{Gtransfo'}
\tilde{x}_i(\textbf{x};\epsilon)=x_i+\epsilon \xi_i ({\textbf{x}})+O(\epsilon ^2), i= 1,\dots,p.
\end{equation}
\end{itemize}
In this paper, we shall assume that the Lie point symmetries are of the form
\begin{equation}\label{Gtransfo''}
\tilde{n}=n;\quad \tilde{u_n} \simeq  u_n+\epsilon Q(n,u_n)
\end{equation}
and that the corresponding infinitesimal generator is given by
\begin{eqnarray}\label{Ngener}
X=  Q(n,u_u)\partial u_n+ SQ(n,u_n)\partial_{u_{ n+1}} +\cdots +S ^{N-1}Q(n,u_n)\partial _{u_{n +N-1}},
\end{eqnarray}
where the shift operator, $S$, is defined as $S:n\mapsto n+1$. \par
\noindent The symmetry condition is given by
\begin{equation}\label{symcdgeneral}
\tilde{u}_{n+N}=\omega(n, \tilde{u}_n, \tilde{u}_{n+1}, \dots, \tilde{u}_{n+N-1})
\end{equation}
whenever (\ref{general}) holds. The substitution of (\ref{Gtransfo''}) into equation (\ref{symcdgeneral}) yields the linearized symmetry condition
\begin{eqnarray}\label{Gsymcdt}
\mathcal{S}^{(N)} Q- X \omega=0
\end{eqnarray}
whenever (\ref{general}) holds. The first integral of (\ref{general}) is given by
\begin{eqnarray}\label{GFintcdt}
(S-id)\phi (n, u_n, u_{n+1}, \dots, u_{n+N-1}) =0
\end{eqnarray}
whenever (\ref{general}) holds.

\noindent It is known that for second-order homogeneous linear difference equations, $u_{n+2}=a(n)u_n + b(n)u_{n+1}$, the Lie algebra of symmetry generator is height dimensional. The characteristics are given by
\begin{align}\label{Qgener}
\begin{split}
Q_1=&  U_1(n),\;  Q_2= U_2(n),\; Q_3= \phi ^1 U_1(n),\; Q_4= \phi ^2 U_1(n),\\ Q_5= &\phi ^1 U_2(n), Q_6= \phi ^2 U_2(n),\; Q_7= (\phi ^1 )^2 U_1(n) \phi ^1 \phi^2 U_2(n)=\phi ^1 u_n,\\ Q_8=& (\phi ^2 )^2 U_2(n)+ \phi ^1 \phi^2 U_1(n)=\phi ^2 u_n,
\end{split}
\end{align}
where $U_1$ and $U_2$ are two linearly independent solutions of the equation, and
\begin{align}
\phi ^1 
=\dfrac{u_n SU_2 - U_2 u_{n+1}}{U_1SU_2-U_2SU_1};\,\quad
\phi ^2
=\dfrac{u_{n+1} U_1 - u_{n}SU_1}{U_1SU_2-U_2SU_1}.
\end{align}
The Lie algebra in this case is isomorphic to $\mathfrak{sl}(3)$. However, these results do not hold for higher order. \par \noindent In this paper, we restrict ourselves to fourth-order difference equations and we introduce a technique for finding their symmetries and conservation laws.
\section{Symmetries}
Consider a fourth-order recurrence equation
\begin{equation}\label{E}
u_{n+4}=\omega(n, u_n, u_{n+1}, \dots, u_{n+3}).
\end{equation}
Then, the linearized symmetry condition (\ref{symcdgeneral}) simplifies to
\begin{eqnarray}\label{symcdt}
\mathcal{S}^{(4)} Q= X \omega,
\end{eqnarray}
with
\begin{eqnarray}\label{gener}
X=  Q(n,u_u)\partial u_n+ SQ(n,u_n)\partial_{u_{ n+1}} +\cdots +S ^3Q(n,u_n)\partial _{u_{n +3}}
\end{eqnarray}
as the corresponding symmetry generator. To find the characteristic we first impose the symmetry condition (\ref{symcdt}) and we then differentiate it with respect to $u_n$ by keeping $\omega$ fixed and by seeing  $u_{n+1}$ as a function of $u_n,u_{n+2},u_{n+3}$ and $\omega$. This leads to
\begin{align}
\begin{split}
&\frac{ \omega_{,u_n}}{ \omega_{,u_{n+1}}}\left[  \omega_{,u_nu_{n+1}}Q-
\left(\omega_{ u_{n+1}} SQ\right)_{,u_{n+1}}-\left(\omega_{,u_{n+1}u_{n+2}}\right)S^2Q-
\left(\omega_{,u_{n+2}u_{n+3}}\right)S^3Q\right]\\&
-\left[  \omega_{,u_n} Q\right]_{,u_n}+\left(\omega_{,u_nu_{n+1}}\right)SQ+\left(\omega_{,u_nu_{n+2}}\right)S^2Q+
\left(\omega_{,u_nu_{n+3}}\right)\;S^3Q=0,
\end{split}
\end{align}
 where $f_{,x}$ denotes the derivative of $f$ with respect to $x$. If $\omega$ is independent of $u_{n+2}$ and $u_{n+3}$, then the above equation simplifies to
\begin{align}\label{deteQ}
\begin{split}
& \big[  \omega_{,u_n} Q\big]_{,u_n}-\left(\omega_{,u_nu_{n+1}}\right)SQ
-\dfrac{ \omega_{,u_n}}{ \omega_{,u_{n+1}}}\Big[ \left(\omega_{,u_nu_{n+1}}\right)Q
-\left(\omega_{ u_{n+1}} SQ\right)_{,u_{n+1}}
\Big]=0.
\end{split}
\end{align}
which we shall refer to as determining equation. In the next two examples we shall use the above determining equation to find symmetries of two fourth-order difference equations.
\subsection{Applications}
Consider the fourth-order difference equation
\begin{equation}\label{symexa1}
u_{n+4}=a(n)u_{n}+b(n)u_{n+1},
\end{equation}
where $a$ and $b$ are functions of $n$.
The reader can readily verify that the solution to the determining equation (\ref{deteQ}) in this case is given by
\begin{align}
Q(n,u_{n})=c_1u_n +c_2(n),
\end{align}
where $c_2$ is a function of $n$ and $c_1$ a constant. The substitution of the characteristic in the symmetry condition (\ref{symcdt}) requires that the function $c_2$ satisfies the original equation (\ref{symexa1}), i.e.,
\begin{align}\label{symexa1'}
c_2(n+4)=a(n)c_2(n)+b(n)c_2(n+1).
\end{align}
If we let $U_1, U_2, U_3$ and $U_4$ be the solutions of (\ref{symexa1'}) then the symmetries are as follows:
\begin{align}
X_0 =& u_n \partial_{u_n},\\
X_1 =& U_1\partial_{u_n},\\ X_2=&U_2\partial _{u_n}, \\ X_3=&U_3 \partial _{u_{n}}, \\ X_4=&U_4 \partial _{u_{n}}.
\end{align}
For instance, if $a(n)=1$ and $b(n)=0$, that is, \begin{align}
u_{n+4}=u_n,
\end{align}
it is easy to check that the solutions to equation (\ref{symexa1'}) in this case are given by $1,\,\cos \frac{n\pi}{2},\,\sin\frac{n\pi}{2}$ and $(-1)^n$. Therefore, the symmetries will be
\begin{equation}\label{25}
\begin{split}
X_{11} =& \partial_{u_n}, X_{12} = (-1)^n \partial_{u_n}, X_{13}=\cos \frac{n\pi}{2}\partial _{u_n},
 X_{14}=\sin\frac{n\pi}{2} \partial _{u_{n}}, \\X_{15} = &u_n\partial_{u_n}.
\end{split}
\end{equation}
We have applied this result to the fourth-order difference equation
\begin{align}\label{exa1'}
u_{n+4}={\dfrac{n}{n+4}}u_n
\end{align}
and we have found that its symmetries are given by
\begin{equation}
\begin{split}
X_0 =& \dfrac{4}{n}\partial_{u_n}, X_1 = \dfrac{4(-1)^n}{n} \partial_{u_n}, X_2=\dfrac{4}{n}\cos \frac{n\pi}{2}\partial _{u_n},  X_3=\dfrac{4}{n}\sin\frac{n\pi}{2} \partial _{u_{n}},\\
 X_4 = &u_n\partial_{u_n}.
\end{split}
\end{equation}
In this next example, we consider an equation obtained from \cite{chen} where the authors looked at the dynamical properties of
$$
u_{n+1}=\frac{u_{n-2}^a+u_{n-3}}{u_{n-2}^au_{n-3}+1}
$$
with positive initial conditions and $a \in [0,1)$. We choose to find symmetries of this fourth-order difference equation when $a=1$, i.e.,
\begin{equation}\label{example2}
u_{n+4}=\frac{u_{n+1}+u_n}{u_nu_{n+1}+1}.
\end{equation}
The determining equation (\ref{deteQ}) becomes
\begin{align}\label{a2}
\begin{split}
& \bigg[\frac{2 \, {\left(u_n + u_{n+1}\right)} }{{\left(u_n u_{n+1} + 1\right)}^{3}} \bigg] SQ- 2 \, \bigg[\frac{{\left({u_{n+1}}^{2}-1\right)} {u_{n+1}}}{{\left(u_nu_{n+1} + 1\right)}^{3}}
 \bigg] Q+ {\left(\frac{{u_{n+1}}^{2}-1}{{\left(u_n u_{n+1} + 1\right)}^{2}} \right)} Q'
+
\\&
 \Bigg[2 \, {\left(\frac{{\left(u_n + u_{n+1}\right)} {u_n}^{2}}{{\left(u_n  u_{n+1} + 1\right)}^{3}} - \frac{u_n}{{\left(u_n  u_{n+1} + 1\right)}^{2}}\right)} SQ
 +
 \bigg(\frac{2 \, {\left(u_n + u_{n+1}\right)} }{{\left(u_n u_{n+1} + 1\right)}^{3}} \bigg) Q
 \\&
 -\bigg(\frac{{u_{n}}^{2}-1}{{\left(u_n u_{n+1} + 1\right)}^{2}}\bigg) (SQ)'\Bigg]. 
\end{split}
\end{align}
After clearing fractions and differentiating three times with respect to $u_n$ we get
\begin{align}\label{a7}
\begin{split}
&{u_n}^{2} { u_{n+1}}^{3}Q^{(4)} + 4 \, u_n { u_{n+1}}^{3} Q^{(3)} - {u_n}^{2} u_{n+1} Q^{(4)}- {u_{n+1}}^{3} Q^{(4)} - \\&4 \, u_n  u_{n+1} Q^{(3)} + u_{n+1} Q^{(4)}=0.
\end{split}
\end{align}
Equation (\ref{a7}) can be split with respect to $u_{n+1}$ as follows
\begin{subequations}\label{splita1}
\begin{eqnarray}
{u_{n+1}}^{3}&:& ({u_n}^{2}-1) Q^{(4)} +4u_nQ^{(3)}=0\\
u_{n+1}&:&-({u_n}^{2}-1) Q^{(4)} -4u_nQ^{(3)}=0
\end{eqnarray}
\end{subequations}
The most general solution of the above system is given by
\begin{align}\label{sol:}
\begin{split}
Q
 =&  c_1(n) \left[  (2{u_n}^2-2)\ln \left(\frac{1-u_n}{1+u_n} \right)-4u_n+1\right]+  c_2+c_3 u_n +c_4 {u_n}^2,
\end{split}
\end{align}
where the $c_i,\; i=1,\dots, 4$, are functions of $n$.
\begin{itemize}
\item
Assuming $c_1=0$, that is, $Q=c_4 {u_n}^2 +c_3 u_n + c_2$, the
substitution of (\ref{sol:}) in (\ref{a2}),(\ref{deteQ}) and (\ref{symcdt}) allows us to find the dependency among the $c_i$. We found that
\begin{align}
c_3(n) =&0, \quad c_4(n)=-c_2(n),
\end{align}
with $c_4$ satisfying the recurrence equation
\begin{align}
c_4(n)+c_4(n+1)-c_4(n+2)=0.
\end{align}
Therefore, the possible characteristics here are $Q_1 =(u_n^2-1)\left( \frac{ 1+\sqrt{5}}{2}\right)^n$ and $Q_2=(u_n^2-1)\left( \frac{ 1-\sqrt{5}}{2}\right)^n$. 
\item
If we assume that $c_1\neq 0$ then the substitution of (\ref{sol:}) in (\ref{a2}),(\ref{deteQ}) and (\ref{symcdt}) leads to
\begin{align}\label{withln}
c_1(n) =&\frac{1}{4}c_3(n),\quad c_2(n)=-\frac{1}{4}c_3(n),\quad c_3(n)=c,\quad c_4(n)=0,
\end{align}
where $c$ is an arbitrary constant. We then substitute (\ref{withln}) in (\ref{sol:}) to  get another characteristic $Q=\frac{1}{2}(u_n^2 -1)\ln |\frac{1-u_n}{1+u_n}|$.
\end{itemize}
In all, (\ref{example2}) has three non trivial symmetries given by
\begin{align}
\begin{split}
X_1 =& (u_n^2-1)\left( \frac{ 1+\sqrt{5}}{2}\right)^n \partial_{u_n},\quad X_2= (u_n^2-1)\left( \frac{ 1-\sqrt{5}}{2}\right)^n\partial _{u_n},\\  X_3= &\frac{1}{2}(u_n^2 -1)\ln \bigg|\frac{1-u_n}{1+u_n}\bigg|\partial _{u_n}.
\end{split}
\end{align}
\section{Conservation Laws}
Let \begin{align}\phi = \phi (n,u_n, u_{n+1}, u_{n+2}, u_{n+3})\end{align} be a first integral of equation (\ref{E}). By definition, $\phi$ is constant on the solution of (\ref{E}). We  then have
\begin{align}\label{SCdH}
\phi\left( n+1, u_{n+1}, u_{n+2}, u_{n+3},\omega\right)=
\phi\left( n, u_n, u_{n+1}, u_{n+2}, u_{n+3}\right).
\end{align}
It can readily be verified that
\begin{align}\label{P_i}
\begin{split}
&P_0 = \left(\frac{\partial \omega}{\partial u_n}\right)\cdot S(P_3),\quad
P_1 = S(P_0)+\left(\frac{\partial \omega}{\partial u_{n+1}}\right)\cdot S(P_3),\quad\\&
P_2 = S(P_1)+\left(\frac{\partial \omega}{\partial u_{n+2}}\right)\cdot S(P_3),\quad
P_3 = S(P_2)+\left(\frac{\partial \omega}{\partial u_{n+3}}\right)\cdot S(P_3),
\end{split}
\end{align}
where
\begin{align}\label{P_33}
P_0=\frac{\partial \phi}{\partial u_n},\quad P_1=\frac{\partial \phi}{\partial u_{n+1}},\quad P_2=\frac{\partial \phi}{\partial u_{n+2}}\quad \text{and} \quad P_3=\frac{\partial \phi}{\partial u_{n+3}}.
\end{align}
To get the determining equation that will allow us to find the first integrals, we combine equations (\ref{P_i}) in a single equation involving $P_3$ only:
\begin{align}\label{DetEqH}
\begin{split}
&S^3\left(\omega_{,u_n}\right)S^4 P_3+
S^2\left(\omega_{,u_{n+1}}\right)S^3 P_3+
S\left(\omega_{,u_{n+2}}\right)S^2 P_3+
\omega_{,u_{n+3}}S P_3 - P_3=0.
\end{split}
\end{align}
Equation (\ref{DetEqH}) can be solved for $P_3$ and hence equation (\ref{P_33}) will enable us to obtain the first integrals. If $\omega$ is independent of $u_{n+2}$ and $u_{n+3}$, then the determining equation in this case is reduced to
\begin{align}\label{DetEqH'}
\begin{split}
&S^3\left(\omega_{,u_n}\right)S^4 P_3+
S^2\left(\omega_{,u_{n+1}}\right)S^3 P_3 - P_3=0.
\end{split}
\end{align}
\subsection{Applications}
Consider the difference equation
\begin{equation}\label{clsexa1}
u_{n+4}=a(n)u_{n}+b(n)u_{n+1}.
\end{equation}
By assuming that $P_3=P_3(n,u_n)$, we found that the solution to the determining equation (\ref{DetEqH'}) in this case is given by  \begin{align}\label{clu4'}P_3(n,u_n)=K(n), \end{align}
where $K$ satisfies the equation
\begin{align}\label{Ki}
K(n)=a(n+3)K(n+4)+b(n+2)K(n+3).
\end{align}
It follows that $P_0=a(n)K_i(n+1)$, $P_1=a(n+1)K_i(n+2)+b(n)K_i(n+1)=K_i(n-2)$, $P_2=a(n+2)K_i(n+3)+b(n+1)K_i(n+2)=K_i(n-1)$ and $P_3=K_i(n)$, where $K_i$, $i=1,\dots 4$ are the solutions of (\ref{Ki}). Therefore, the first integrals of equations of the form (\ref{clsexa1}) are giving by
\begin{align}\label{clsexa1g}
\begin{split}
\phi_i=&\int\bigg( a(n)K_i(n+1) d u_n +K_i(n-2) d u_{n+1} +  K_i(n-1) d u_{n+2}\\&+ P_3 d u_{n+3}\bigg)+ G_i(n),\,  i=1,\dots 4,
\end{split}
\end{align}
for some function $G_i$.\par \noindent
For the sake of clarification, let us consider the case where $a(n)=1$ and $b(n)=0$
, that is,
\begin{align}\label{47}
u_{n+4}=u_n.
\end{align}
The possible solutions to (\ref{Ki}) are $1,(-1)^n,\cos{\frac{n\pi}{2}}$ and $\sin{\frac{n\pi}{2}}$.
\begin{itemize}
\item If $P_3(n,u_n)=1$ then $P_2(n,u_n)=P_1(n,u_n)=P_0(n,u_n)=1$. So $$\phi (n,u_n, u_{n+1},u_{n+2},u_{n+3})=u_n+u_{n+1}+u_{n+2}+u_{n+3}+c_1(n)$$ for some function $c_1$ of $n$. Imposing $S \phi =\phi$, we found that the conservation law in this case is given by
  \begin{align}
  \phi =u_n+u_{n+1}+u_{n+2}+u_{n+3}+c_1.
  \end{align}
  \item Similarly, if $P_3 =(-1)^n$ then $P_2=(-1)^{n+1}$, $P_1=(-1)^{n}$ and $P_0 = (-1)^{n+1}$.
\item
If $P_3=\cos{\frac{n\pi}{2}}$ then $P_2=\sin{\frac{n\pi}{2}}$, $P_1=-\cos{\frac{n\pi}{2}}$ and $P_0 = -\sin{\frac{n\pi}{2}}$.
\item
If $P_3=\sin{\frac{n\pi}{2}}$ then $P_2=-\cos{\frac{n\pi}{2}}$, $P_1=-\sin{\frac{n\pi}{2}}$ and $P_0=\cos{\frac{n\pi}{2}}$.
\end{itemize}
We therefore obtained four conservation laws for the equation $u_{n+4}=u_n$ given by 
\begin{subequations}\label{clu44}
\begin{align}
  \phi_1 =&u_n+u_{n+1}+u_{n+2}+u_{n+3}\\
  \phi_2=&(-1)^n\left( -u_n+u_{n+1}-u_{n+2}+u_{n+3}\right)\\ 
  \phi_3=&\sin{\frac{n\pi}{2}}\left( u_{n+2}-u_{n}\right)+\cos{\frac{n\pi}{2}}\left( u_{n+3}-u_{n+1}\right)\\
  \phi_4=&\sin{\frac{n\pi}{2}}\left( u_{n+3}-u_{n+1}\right)-\cos{\frac{n\pi}{2}}\left( u_{n+2}-u_{n}\right).
  \end{align}
\end{subequations}
Again, we consider equation (\ref{exa1'}), i.e.,
\begin{align}\label{Exa1'}
u_{n+4}={\dfrac{n}{n+4}}u_n.
\end{align}
Here, condition (\ref{Ki}) becomes
\begin{equation}\label{exa1''}
K_{n+4}=\frac{n+7}{n+3}\,K_n.
\end{equation}
We have proved that $\dfrac{1}{3}(n+3)$, $\dfrac{(-1)^n}{3}(n+3)$, $\dfrac{1}{3}(n+3)\sin\frac{n\pi}{2}$ and $\dfrac{1}{3}(n+3)\cos\frac{n\pi}{2}$ are solutions of (\ref{exa1''}). Thus, the first integrals of (\ref{Exa1'}) are as follows:
\begin{subequations}\label{clu441}
\begin{align}
\begin{split}
  \phi_1 =&\dfrac{1}{3}\bigg[nu_n+(n+1)u_{n+1}+(n+2)u_{n+2}+(n+3)u_{n+3}\bigg]\end{split}\\
  \begin{split}
  \phi_2 =&\dfrac{(-1)^n}{3}\bigg[-nu_n+(n+1)u_{n+1}-(n+2)u_{n+2}+(n+3)u_{n+3}\bigg]
  \end{split}\\
  \begin{split}
  \phi_3 =&\dfrac{1}{3}\bigg[-n\sin{\bigg(\dfrac{n\pi}{2}\bigg)}u_n-(n+1)\cos{\bigg(\dfrac{n\pi}{2}\bigg)}u_{n+1}
  \\&+(n+2)\sin{\bigg(\dfrac{n\pi}{2}\bigg)}u_{n+2}+(n+3)\cos{\bigg(\dfrac{n\pi}{2}\bigg)}u_{n+3}\bigg]
  \end{split}\\
  \begin{split}
  \phi_4 =&\dfrac{1}{3}\bigg[n\cos{\bigg(\dfrac{n\pi}{2}\bigg)}u_n-(n+1)\sin{\bigg(\dfrac{n\pi}{2}\bigg)}u_{n+1}
  \\&-(n+2)\cos{\bigg(\dfrac{n\pi}{2}\bigg)}u_{n+2}+(n+3)\sin{\bigg(\dfrac{n\pi}{2}\bigg)}u_{n+3}\bigg].
  \end{split}
  \end{align}
\end{subequations}
\par  \noindent In the previous examples, we supposed that $P_3=P_3(n,u_n)$, i.e., we assumed that
\begin{align}\label{fg}\phi = f(n,u_n)u_{n+3}+g(n,u_n,u_{n+1},u_{n+2}).
\end{align}
It has to be noted that it is not every fourth-order difference equation that has first integrals of this form. For example, equation (\ref{example2}) given by
\begin{equation}
u_{n+4}=\dfrac{u_{n+1}+u_n}{u_nu_{n+1}+1}
\end{equation}
does not have first integrals of the form (\ref{fg}). We then decided to investigate the existence of first integrals of (\ref{example2}) where $\dfrac{\partial\phi}{\partial u_{n+3}}=P_3(n,u_{n+1})$. With this assumption, the determining equation (\ref{DetEqH'}) reduces to
\begin{align}\label{DetEqexa2}
\begin{split}
&S^3\left( \omega_{ u_n}\right)S^4 P_3(n+1,u_{n+1})+
S^2\left( \omega_{ u_{n+1}}\right)S^3 P_3(n+1,u_{n+1}) \\&- P_3(n+1,u_{n+1})=0.
\end{split}
\end{align}
The operator $S^{-1}:n\mapsto n-1$ acts on (\ref{DetEqexa2}) to produce
\begin{align}\label{DetEqexa2'}
\begin{split}
&S^2\left( \omega_{ u_n}\right) S^4P_3(n,u_n)+
S\left( \omega_{ u_{n+1}}\right)S^3P_3(n,u_{n}) - P_3(n,u_{n})=0.
\end{split}
\end{align}
In the above equation, the function $P_3$ takes different arguments. To solve this, we first differentiate it with respect to $u_n$, keeping $\omega$ fixed, and by assuming that $u_{n+3}$ is a function of $u_n$ and $\omega$.
Secondly, we differentiate the resulting equation with respect to $u_{n}$ to get
\begin{align}
\begin{split}
&\big[{u_{n+3}}^{2} {u_n}^{2} { u_{n+1}}^{3} {u_{n+2}}^{5}+3 \, {u_{n+3}}^{2} {u_n}^{2} { u_{n+1}}^{2} {u_{n+2}}^{4} + 2 \, {u_{n+3}} {u_n}^{2} { u_{n+1}}^{3} {u_{n+2}}^{4} \\
 &- {u_{n+3}}^{2} { u_{n+1}}^{3} {u_{n+2}}^{5}  + +3 \, {u_{n+3}}^{2} {u_n}^{2} { u_{n+1}} {u_{n+2}}^{3}  +6 \, {u_{n+3}} {u_n}^{2} { u_{n+1}}^{2} {u_{n+2}}^{3} \\
  &+ {u_n}^{2} { u_{n+1}}^{3} {u_{n+2}}^{3}  - 3 \, {u_{n+3}}^{2} { u_{n+1}}^{2} {u_{n+2}}^{4}  - 2 \, {u_{n+3}} { u_{n+1}}^{3} {u_{n+2}}^{4}  \\
 &  + {u_{n+3}}^{2} {u_n}^{2} {u_{n+2}}^{2}  + 6 \, {u_{n+3}} {u_n}^{2} { u_{n+1}} {u_{n+2}}^{2}  + 3 \, {u_n}^{2} { u_{n+1}}^{2} {u_{n+2}}^{2}  - \\
    &3 \, {u_{n+3}}^{2} { u_{n+1}} {u_{n+2}}^{3}  - 6 \, {u_{n+3}} { u_{n+1}}^{2} {u_{n+2}}^{3}  - { u_{n+1}}^{3} {u_{n+2}}^{3}  \\
    &+2 \, {u_{n+3}} {u_n}^{2} {u_{n+2}}  + 3 \, {u_n}^{2} { u_{n+1}} {u_{n+2}} -{u_{n+3}}^{2} {u_{n+2}}^{2}  - 6 \, {u_{n+3}} { u_{n+1}} {u_{n+2}}^{2} -\\
      &3 \, { u_{n+1}}^{2} {u_{n+2}}^{2}   + {u_n}^{2}  - 2 \, {u_{n+3}} {u_{n+2}} - 3 \, { u_{n+1}} {u_{n+2}}  - 1
 +\big]P_3''(n,u_{n}) \\
 &+ \big[2 \, {u_{n+3}}^{2} {u_n} { u_{n+1}}^{3} {u_{n+2}}^{5}+ 6 \, {u_{n+3}}^{2} {u_n} { u_{n+1}}^{2} {u_{n+2}}^{4}  + \\
 &4 \, {u_{n+3}} {u_n} { u_{n+1}}^{3} {u_{n+2}}^{4}   +6 \, {u_{n+3}}^{2} {u_n} { u_{n+1}} {u_{n+2}}^{3} + 12 \, {u_{n+3}} {u_n} { u_{n+1}}^{2} {u_{n+2}}^{3} +\\
    &2 \, {u_n} { u_{n+1}}^{3} {u_{n+2}}^{3}+ 2 \, {u_{n+3}}^{2} {u_n} {u_{n+2}}^{2} + 12 \, {u_{n+3}} {u_n} { u_{n+1}} {u_{n+2}}^{2}  +\\
     &6 \, {u_n} { u_{n+1}}^{2} {u_{n+2}}^{2} + 4 \,{u_{n+3}} {u_n} {u_{n+2}} + 6 \, {u_n} { u_{n+1}} {u_{n+2}}    + 2 \, {u_n}      \big] P_3'(n,u_{n}).
\end{split}
\end{align}
We can now split with respect to multiples of powers of $u_{n+1}$,$u_{n+2}$ and $u_{n+2}$. This leads a single differential equation
\begin{equation}
({u_n}^2-1)P_3''(n,u_n)+2u_nP_3'(n,u_n)=0
\end{equation}
whose solution is $P_3(n,u_n)=c_1\log\sqrt{\dfrac{1-u_n}{1+u_n}}+c_2$ and therefore we can assume without loss of generality that
\begin{align}
P_3(n,u_{n+1})=&c_1(n)\log\sqrt{\dfrac{1-u_{n+1}}{1+u_{n+1}}}+c_2(n).
\end{align}
\section{Multipliers and Associations}
In this section we look at the more recent notion of multipliers (also known as characteristics) whose role is to identify equivalent conservation laws and we check whether or not the symmetries and first integrals that we obtained before are 'associated'. 
%
We use the shift operator and the first integrals to obtain the multiplier. This method is more direct but works for a certain class of difference equations only. A standard way for finding multipliers can be found in \cite{12}.\par \noindent Recall that the first integrals of (\ref{47}) are given by (\ref{clu44}). We have verified that
\begin{subequations}
\begin{align}
(S-id)\Phi_{11} = &[1](E_1)\\
(S-id)\Phi_{12} = &[(-1)^{n+1}](E_1)                            \\
(S-id)\Phi_{13} = &\left[-\sin{\left(\dfrac{n\pi}{2}\right)}\right](E_1)\\
(S-id)\Phi_{14} = &\left[\cos{\left(\dfrac{n\pi}{2}\right)}\right](E_1).
\end{align}
\end{subequations}
We conclude that $(E_1)$ has four multipliers given by $1$, $(-1)^{n+1}$, $-\sin{\left(\dfrac{n\pi}{2}\right)}$ and $\cos{\left(\dfrac{n\pi}{2}\right)}$. These multipliers turn out to be the four characteristics ( of the symmetries) of $(E_1)$, given in (\ref{25}), we found earlier. 
\par \noindent
Similarly, we have verified that for ($E_2$) we have
\begin{subequations}
\begin{align}
(S-id)\Phi_{21} = &\left[\frac{n+4}{3}\right](E_2)\\
(S-id)\Phi_{22} = &\left[\dfrac{(-1)^{n+1}}{3}(n+4)\right](E_2)                            \\
(S-id)\Phi_{23} = &\left[-\frac{n+4}{3}\sin{\left(\dfrac{n\pi}{2}\right)}\right](E_2)\\
(S-id)\Phi_{24} = &\left[\frac{n+4}{3}\cos{\left(\dfrac{n\pi}{2}\right)}\right](E_2).
\end{align}
\end{subequations}
Thus the four multipliers of equation ($E_2$) are given by $\frac{n+4}{3}$, $\dfrac{(-1)^{n+1}}{3}(n+4)$, $-\frac{n+4}{3}\sin{\left(\dfrac{n\pi}{2}\right)}$ and $\frac{n+4}{3}\cos{\left(\dfrac{n\pi}{2}\right)}$. This implies that the first integrals $\Phi_{1i}$, $1\leq i\leq 4$,(resp. $\Phi_{1i}$) are non-equivalent since they generated different multipliers.\par \noindent
On the other hand, we check the association of the symmetries and first integrals. It is known that when these two tools are associated one could use them to perform double reductions of the equations. The results below show the associated and non-associated symmetries and conserved vectors.
\begin{align}
\begin{split}
& X_{11} \Phi_{11}=4,\quad X_{11} \Phi_{12}=0,\quad X_{11} \Phi_{13}=0,\quad X_{11} \Phi_{14}=0,\\
& X_{12} \Phi_{11}=0,\quad X_{12} \Phi_{12}=-4,\quad X_{12} \Phi_{13}=0,\quad X_{12} \Phi_{14}=0,\\
& X_{13} \Phi_{11}=0,\quad X_{13} \Phi_{12}=0,\quad X_{13} \Phi_{13}=0,\quad X_{13} \Phi_{14}=2,\\
& X_{14} \Phi_{11}=0,\quad X_{14} \Phi_{12}=0,\quad X_{14} \Phi_{13}=-2,\quad X_{14} \Phi_{14}=0,\\
& X_{15} \Phi_{1i}=\Phi_{1l},\quad l=1,\dots,4,
\end{split}
\end{align}
and
\begin{align}
\begin{split}
& X_{21} \Phi_{21}=\dfrac{16}{3},\quad X_{21} \Phi_{22}=0,\quad X_{21} \Phi_{23}=0,\quad X_{21} \Phi_{24}=0,\\
& X_{22} \Phi_{21}=0,\quad X_{22} \Phi_{22}=-\dfrac{16}{3},\quad X_{22} \Phi_{23}=0,\quad X_{22} \Phi_{24}=0,\\
& X_{23} \Phi_{21}=0,\quad X_{23} \Phi_{22}=0,\quad X_{23} \Phi_{23}=0,\quad X_{23} \Phi_{24}=\dfrac{8}{3},\\
& X_{24} \Phi_{21}=0,\quad X_{24} \Phi_{22}=0,\quad X_{24} \Phi_{23}=-\dfrac{8}{3},\quad X_{24} \Phi_{24}=0,\\
& X_{25} \Phi_{2i}=\Phi_{2l},\quad l=1,\dots,4.
\end{split}
\end{align}


\section{Conclusion}
We have presented a technique for obtaining symmetries and conservation laws of fourth-order difference equations. To ease our computation we have made some assumptions that enabled us to derive a number of symmetries and first integrals. 
We completed our investigation by looking at the association of symmetries and first integrals that we have found and multipliers were derived.

\end{document}